\documentclass[12pt]{article}

\newcommand{\nc}{\newcommand}
\nc{\dfrac}{\displaystyle \frac } 
\newcommand{\be}{\begin{equation}}
\newcommand{\ee}{\end{equation}}
\newcommand{\br}{\begin{array}}
\newcommand{\er}{\end{array}}

\input{tcilatex}

\begin{document}

\begin{center}
{\Large \textbf{Generalized Exclusion and Hopf Algebras}}
\end{center}

\vspace{5mm}

\begin{center}
A. Yildiz

\noindent Istanbul Technical University, Faculty of Science and Letters,
Physics Department, 80626, Maslak, Istanbul, Turkey \footnote{%
E-mail : yildiz@gursey.gov.tr}. \\[0pt]
\end{center}

\vspace{5mm}

\begin{center}
\textbf{Abstract}
\end{center}

We propose a generalized oscillator algebra at roots of unity with
generalized exclusion and investigate the braided Hopf structure. We find
that there are two solutions one of which is the generalized exclusion of
the bosonic type and the other is the generalized exclusion of the fermionic
type. We also discuss the covariance properties of these oscillators.

\bigskip

\noindent

\section{Introduction}

The generalization of the statistics of identical particles can be
done in two different ways in principle. The first approach is
based on the generalization of the $\pm 1$ factor of boson/fermion
algebra which comes from the interchange of two identical
particles. This approach lead to the introduction of anyons
\cite{art1} which found applications in fractional quantum Hall
effect and superconductivity \cite {art2}. The quon algebras
\cite{art3} and braided oscillators \cite{art4} are other
generalizations based on the generalizations of the exchange
parameter. Mathematically this corresponds to the generalization
of the algebra in tensor product space such that

\begin{equation}  \label{eq:perm}
(1\otimes a)(b\otimes 1)=\pm \pi (a\otimes b)=\pm b\otimes a
\end{equation}

\noindent (where the + sign refers to bosons and - sign refers to fermions)
the permutation map $\pm \pi $\ is replaced by a generalized map $\psi $ i.e,

\begin{equation}  \label{eq:bra}
(1\otimes a)(b\otimes 1)=\psi (a\otimes b).
\end{equation}

\noindent In general we have $\psi ^{2}\neq id$ and hence the map $\psi $ is
called braiding. The algebra in the tensor product space is connected with
Hopf algebras such that the coproduct map $\Delta $ of Hopf algebras is a
homomorphism to the tensor product space $\Delta :A\rightarrow A\otimes A.$
It creates identical copies of algebras and hence the theory of Hopf
algebras can be used in the operator algebra formulations of identical
particles. For bosons/fermions ordinary/super Hopf algebras are used and for
the generalized exchange phase braided Hopf algebra axioms \cite{art5} are
used. The braided Hopf algebra axioms are generalizations of ordinary Hopf
algebra axioms such that in the limit $\psi \rightarrow \pi $ braided Hopf
algebra axioms reduce to the ordinary Hopf algebra axioms.

The other approach for the generalization of statistics is based on the
generalization of the Pauli exclusion principle put forward by Haldane \cite
{art6} and the statistical mechanics is developed by Wu \cite{art7}. In the
literature an extensive amount of work has been done as applications of this
generalization. Our aim in this paper is to find an operator algebra for
this generalized exclusion. It is known that the $q$-deformed generalization
of nilpotent algebras requires the deformation parameter to be a root of
unity \cite{art12}. For example the q-deformed harmonic oscillators at roots
unity is connected with two anyon system \cite{art8}. Different versions of
the deformation of oscillator algebras and their use in physics can be found
in the review \cite{art9} and the connection of between q-deformed
structures and fractional statistics can be found in \cite{art10}.\ In terms
of creation and annihilation operators the generalized exclusion can be
expressed as $a^{k}=0$ and $(a^{\dagger })^{k}=0$. The limit $k\rightarrow 2$
gives the Pauli exclusion principle for fermions. Some operator algebras for
fractional exclusion statistics are proposed in ref. \cite{art11}.

In this work we discuss the problem of the construction of an oscillator
algebra unifying the notions of generalized phase (braiding), Hopf algebra
and exclusion statistics, i.e., we investigate the braided Hopf algebra
structure of a generalized oscillator with generalized exclusion $a^{k}=0.$

\section{Generalized exclusion algebra}

We start with a generalized oscillator algebra generated by ($a,a^{\ast
},q^{N},q^{-N},1$) satisfying
\begin{eqnarray}  \label{eq:osc1}
aa^{\ast }-Q_{1}a^{\ast }a &=&Q_{2}q^{2N}+Q_{3}q^{-2N},  \nonumber \\
aa^{\ast }-Q_{1}^{-1}a^{\ast }a &=&Q_{2}^{\ast }q^{-2N}+Q_{3}^{\ast }q^{2N},
\nonumber \\
aq^{N} &=&qq^{N}a,\quad q^{N}a^{\ast }=qa^{\ast }q^{N}, \\
aq^{-N} &=&q^{-1}q^{-N}a,\quad q^{-N}a^{\ast }=q^{-1}a^{\ast }q^{-N}
\nonumber \\
q^{N}q^{-N} &=&q^{-N}q^{N}=1_{A},\quad a^{k}=(a^{\ast })^{k}=0.
\nonumber
\end{eqnarray}

\noindent where $q$ and $Q_{1}$ are complex parameters of unit modulus and $%
Q_{2}$ and $Q_{3}$ are any complex parameters and $1_{A}$ is the
identity of the algebra. For the *-structure we impose

\be\label{eq:star2}
 (a^{\ast })^{\ast }=a,\quad  \ (q^{\pm
N})^{\ast }=q^{\mp N}. \ee

\noindent  From the  relations (\ref{eq:osc1}) we get
\begin{eqnarray}  \label{eq:osct}
a^{\ast }a&=&\dfrac{(Q_{2}-Q_{3}^{\ast })q^{2N}+(Q_{3}-Q_{2}^{\ast })q^{-2N}%
}{Q_{1}^{-1}-Q_{1}}  \nonumber \\
aa^{\ast }&=&\dfrac{(Q_{1}^{-1}Q_{2}-Q_{1}Q_{3}^{\ast
})q^{2N}+(Q_{1}^{-1}Q_{3}-Q_{1}Q_{2}^{\ast })q^{-2N}}{Q_{1}^{-1}-Q_{1}}.
\end{eqnarray}

\noindent We note that in contrast to the operator algebra
proposed by Karabali-Nair \cite{art11} where the number operator
is function of the raising and lowering operators,i.e.,
$N=f(a^{\ast }a)$, we take the Hermitian operator $a^{\ast }a)$ as
a function of the  number operator, i.e., $a^{\ast }a=f(N)$. Since
$a$ and $a^{\ast }$ are lowering and raising operators
respectively the operator $ aa^{\ast }$ should satisfy $ aa^{\ast
}=f(N+1)$. The function $f(N)$ for the algebra we propose is given
by (\ref{eq:osct}).

\begin{equation}
q^{2}(Q_{2}-Q_{3}^{\ast })=Q_{1}^{-1}Q_{2}-Q_{1}Q_{3}^{\ast }.
\end{equation}

\noindent The special cases $Q_{1}=q^{-2}$ \ and $Q_{1}=q^{2}$, with a
rescaling of the generators, reduce to

\begin{equation}
aa^{\ast }-q^{-2}a^{\ast }a=q^{2N},\ \ \ aa^{\ast }-q^{2}a^{\ast }a=q^{-2N}
\end{equation}

\noindent which is the well known q-oscillator algebra \cite{art13}. At
roots of unity this algebra is related to two anyon system \cite{art8}. For $%
Q_{1}\neq q^{2},q^{-2}$ we express $Q_{3}$ in terms of $Q_{2}$. Then the
oscillator algebra turns out to be
\begin{eqnarray}
aa^{\ast }-Q_{1}a^{\ast }a=Q_{2}q^{2N}+Q_{1}RQ_{2}^{\ast }q^{-2N}  \nonumber
\\
aq^{N}=qq^{N}a,\ \ q^{N}a^{\ast }=qa^{\ast }q^{N},\ \ a^{k}=0
\end{eqnarray}
and their *-conjugates. The parameter $R\equiv \dfrac{1-q^{2}Q_{1}}{%
Q_{1}-q^{2}}$\bigskip is real. For the Fock space representation we take
\begin{equation}
a\mid n\rangle =a_{n}\mid n-1\rangle ,\quad a^{\ast }\mid n\rangle
=a_{n+1}^{\ast }\mid n+1\rangle ,\quad q^{N}\mid n\rangle =q^{n}\mid
n\rangle .
\end{equation}

\noindent Solving the difference equations we get the spectrum
\begin{equation}
a^{\ast }a\mid n\rangle =\dfrac{Q_{2}q^{2n}-Q_{2}^{\ast }Q_{1}q^{-2(n-1)}}{%
q^{2}-Q_{1}}\mid n\rangle .
\end{equation}

\noindent The existence of a ground state $a\mid 0\rangle =0$ implies
\begin{equation}  \label{eq:par1}
Q_{2}^{\ast }=q^{-2}Q_{1}^{-1}Q_{2}.
\end{equation}

\noindent Since $q$ and $Q_{1}$ are complex parameters with unit modulus we
write the parameters as

\begin{equation}
Q_{1}=e^{i\alpha },\quad q=e^{i\beta },\quad Q_{2}=Be^{i\gamma }
\end{equation}

\noindent and substitute these into (\ref{eq:par1}) to get $e^{i2\gamma
}=e^{i(\alpha +2\beta )}.$ Substituting $Q_{2}$ into the oscillator algebra
we get

\begin{eqnarray}
aa^{\ast }-Q_{1}a^{\ast }a=BqQ_{1}{}^{1/2}q^{2N}+Bq^{-1}Q_{1}^{1/2}Rq^{-2N}
\nonumber \\
aa^{\ast }-Q_{1}^{-1}a^{\ast
}a=Bq^{-1}Q_{1}{}^{-1/2}q^{-2N}+BqQ_{1}^{-1/2}Rq^{2N}
\end{eqnarray}

\noindent Then the spectrum turns out to be

\begin{equation}
a_{n}^{\ast }a_{n}=B\dfrac{\sin (2n\beta )}{\sin \QDOVERD( ) {2\beta -\alpha
}{2}}=B\dfrac{q^{2n}-q^{-2n}}{qQ_{1}^{-1/2}-q^{-1}Q_{1}^{1/2}}.
\end{equation}

\noindent The positive definiteness of $a_{n}^{\ast }a_{n}$ implies that the
representation should be finite dimensional with the states$\mid 0\rangle
......\mid k-1\rangle $ and the highest state should be annihilated by the
raising operator $a^{\ast }\mid k-1\rangle =0.$ This gives that $q$ is a
root of unity $q=e^{i\frac{\pi }{2k}}$ and taking $\dfrac{\pi }{k}-2\pi
\langle \alpha \langle \dfrac{\pi }{k}$ the spectrum turns out to be
\begin{equation}
a_{n}^{\ast }a_{n}=B\dfrac{\sin (n\dfrac{\pi }{k})}{\sin \frac{1}{2}(\dfrac{%
\pi }{k}-\alpha )}.
\end{equation}

\noindent Although the deformation parameter $Q_{1}$ (or $\alpha $) seems
only as a free scaling parameter in the spectrum, its value is going to be
fixed by the Hopf algebra structure.

\section{Braided Hopf Algebra}

The generalization of the permutation map of boson/fermion algebras in
tensor product space by a generalized map leads naturally to the
generalization of the Hopf algebra called braided Hopf algebra whose axioms
in algebraic (not diagrammatic) form read as

\begin{eqnarray}  \label{eq:bhopf}
m\circ (id\otimes m)&=&m\circ (m\otimes id)  \nonumber \\
m\circ (id\otimes \eta )&=&m\circ (\eta \otimes id) = id  \nonumber \\
(id\otimes \Delta )\circ \Delta &=&(\Delta \otimes id)\circ \Delta  \nonumber
\\
(\epsilon \otimes id)\circ \Delta &=&(id\otimes \epsilon )\circ \Delta = id
\nonumber \\
m\circ (id\otimes S)\circ \Delta &=&m\circ (S\otimes id)\circ \Delta =\eta
\circ \epsilon  \nonumber \\
\psi \circ (m \otimes id) &=& (id \otimes m)\circ (\psi \otimes id)\circ (id
\otimes \psi)  \nonumber \\
\psi \circ (id \otimes m) &=& (m \otimes id)\circ (id \otimes \psi)\circ
(\psi \otimes id)  \nonumber \\
(id \otimes \Delta)\circ \psi &=& (\psi \otimes id)\circ (id \otimes \psi
)\circ (\Delta \otimes id)  \nonumber \\
(\Delta \otimes id)\circ \psi &=& (id \otimes \psi )(\psi \otimes id)\circ
(id \otimes \Delta)  \nonumber \\
\Delta \circ m &=& (m \otimes m)(id\otimes \psi \otimes id)\circ (\Delta
\otimes \Delta ) \\
S\circ m &=& m\circ \psi \circ (S\otimes S)  \nonumber \\
\Delta \circ S &=& (S\otimes S)\circ \psi \circ \Delta  \nonumber \\
\epsilon \circ m &=& \epsilon \otimes \epsilon  \nonumber \\
(\psi \otimes id)\circ (id \otimes \psi)\circ (\psi \otimes id)&=&(id
\otimes \psi)\circ (\psi \otimes id)\circ (id \otimes \psi)  \nonumber
\end{eqnarray}

\noindent where $m:A\otimes A\rightarrow A$ is the multiplication map, $%
\Delta :A\rightarrow A\otimes A$ is the comultiplication map,
$\eta :K\rightarrow A$ is the unit map, $\epsilon :A\rightarrow K$
is the counit map, $S:A\rightarrow A$ is the antipode map and
$\psi :A\otimes A\rightarrow A\otimes A$ is the braiding map. Note
that in the limit $\psi \rightarrow \pm \pi $ the braided Hopf
algebra axioms reduce to the ordinary/super Hopf algebra axioms.
The consistency of the braided Hopf algebra axioms
(\ref{eq:bhopf}) requires that the identity element $1_{A}$ in any
algebra should satisfy the following conditions

\be \label{eq:iden}
 \Delta (1_{A})=1_{A}\otimes 1_{A},\ S(1_{A})=1_{A},\ \epsilon(1_{A})=1,\
\psi (1_{A}\otimes a)=a\otimes 1_{A},\ \psi (a\otimes
1_{A})=1_{A}\otimes a \ \ \forall a\in A. \ee

\noindent The $\ast $-structure for the braided algebra $B$
satisfies

\begin{eqnarray}
\Delta \circ \ast &=&\pi \circ (\ast \otimes \ast )\circ \Delta  \nonumber \\
S\circ \ast &=&\ast \circ S.
\end{eqnarray}

\noindent If we interpret the Hopf algebra to define a system (of
oscillators e.g.) then the *-structure should be defined in such a way that
it should be compatible with the algebra in the tensor product space. If $%
a_{1}\equiv a\otimes 1_{A}$ and $b_{2}\equiv 1_{A}\otimes b$ then

\begin{equation} \label{eq:star}
(a\otimes b)^{\ast }=((a\otimes 1)(1\otimes b))^{\ast }=(a_{1}b_{2})^{\ast
}=b_{2}^{\ast }a_{1}^{\ast }=\psi (b^{\ast }\otimes a^{\ast })
\end{equation}

\noindent where in the non-braided limit it gives the ordinary *-structure $%
(a\otimes b)^{\ast }=a^{\ast }\otimes b^{\ast }$. We also note
that the *-involution we define  by (\ref{eq:star}) for the tensor
product is different from the *-involution defined by Majid
\cite{art5}.

To find the Hopf algebra structure of the oscillator algebra
(\ref{eq:osc1}) we start with the general forms of the coproducts
and braidings which are linear in both of the factors of the
tensor product. Although the nonlinear factors are also possible
theoretically we do not consider them in this work. Thus we start
with the most general linear forms of the coproducts, braidings
and antipodes, e.g.,

\be\label{eq:hopf7} \Delta (q^N) = D_1q^N\otimes q^N
+F_1q^N\otimes q^{-N}+G_1a\otimes a^{*} +H_1a^{*}\otimes
a+K_11_{A}\otimes q^N ...
 \ee
\be \psi (q^N \otimes q^N) = g_1q^N\otimes q^N+ g_2a\otimes a^{*}
+g_3a^{*}\otimes a+g_41_{A}\otimes q^N+g_61_{A}\otimes1_{A}+...
\ee

 \noindent and similar expressions for the other generators of the algebra.
It seems that it is very hard to solve the equations for the
braided Hopf algebra because there are too many parameters. When
we substitute the general forms into the relation
$q^Nq^{-N}=q^{-N}q^N=1_{A}$ and use (\ref{eq:iden}) we obtain that
the braidings of $q^N$ and $q^{-N}$ with other generators are
trivial, i.e. they are bosonic elements. From these relations we
also obtain that

\be\label{eq:number1} \Delta (q^{N})=D_{1}q^{N}\otimes q^{N},\quad
\Delta (q^{-N})=D_{1}^{-1}q^{-N}\otimes q^{-N},\quad
S(q^{N})=S_{1}q^{-N},\quad S(q^{-N})=S_{1}^{-1}q^{N} \ee

\noindent where $D_{1}$ and $S_{1}$ are complex parameters. This
simple form (\ref{eq:number1}) simplifies the solution a lot.
Using other relations of the algebra we obtain that the
coproducts, the braidings and the antipodes   for $a$ and $a^{*}$
are in the form

\begin{eqnarray}  \label{eq:hopf3}
 \Delta (a)=D_{2}a\otimes
q^{N}+D_{3}q^{-N}\otimes a,\ \
 \Delta (a^{\ast
})=D_{3}^{\ast }a^{\ast }\otimes q^{N}+D_{2}^{\ast }q^{-N}\otimes
a^{\ast }, \nonumber \\
 S(a)=-S_{2}a,  \quad S(a^{\ast })=-S_{2}^{\ast }a^{\ast },\\
 \psi
(a\otimes a)=za\otimes a,\ \ \psi (a^{\ast }\otimes a^{\ast })=
za^{\ast }\otimes a^{\ast }, \nonumber \\
\psi (a\otimes a^{\ast
})=z^{-1}a^{\ast }\otimes a, \ \  \psi (a^{\ast }\otimes
a)=z^{-1}a\otimes a^{\ast }  \nonumber
\end{eqnarray}

\noindent where $D_{2}$, $D_{3}$, $S_{2}$ are complex parameters
and $z$ is a complex parameter of unit modulus.

\noindent We substitute these  forms into the braided Hopf algebra
axioms (\ref{eq:bhopf}) and find that there are only two solutions
one of which is the bosonic generalization of the exclusion with
the exchange phase $z=Q_{1}=1$ and the other is the fermionic
generalization of the exclusion with the exchange phase
$z=Q_{1}=-1$.

\subsection{Generalized exclusion of the bosonic type}

For the first solution we get $Q_{1}=1$, and all braidings turn out to be
trivial ($z=1$). The oscillator algebra
\begin{eqnarray}  \label{eq:bos1}
aa^{\ast }-a^{\ast }a=Q_{2}(q^{2N}+q^{-2(N+1)}),  \nonumber \\
aq^{N}=qq^{N}a,\quad q^{N}a^{\ast }=qa^{\ast }q^{N},  \nonumber \\
aq^{-N}=q^{-1}q^{-N}a,\quad q^{-N}a^{\ast }=q^{-1}a^{\ast }q^{-N} \\
q^{N}q^{-N}=q^{-N}q^{N}=1_{A},\quad a^{k}=(a^{\ast })^{k}=0
\nonumber
\end{eqnarray}

\noindent has the Hopf algebra structure
\begin{eqnarray}
\Delta (q^{N})=D_{1}q^{N}\otimes q^{N},\ \Delta
(q^{-N})=D_{1}^{-1}q^{-N}\otimes q^{-N},  \nonumber \\
\Delta (a)=D_{2}a\otimes q^{N}+D_{1}^{-2}D_{2}q^{-N}\otimes a,  \nonumber \\
\Delta (a^{\ast })=D_{1}^{2}D_{2}^{\ast }a^{\ast }\otimes q^{N}+D_{2}^{\ast
}q^{-N}\otimes a^{\ast }, \\
S(q^{N})=D_{1}^{-2}q^{-N},\ S(q^{-N})=D_{1}^{2}q^{N},  \nonumber \\
S(a)=-q^{-1}a,\ S(a^{\ast })=-qa^{\ast }  \nonumber \\
\epsilon (q^{N})=D_{1}^{-1},\ \epsilon (q^{-N})=D_{1},\ \epsilon
(a)=\epsilon (a^{\ast })=0  \nonumber
\end{eqnarray}

\noindent where $q=e^{i\frac{\pi }{2k}}$, $D_{1}^{4}=-q^{2}$, $D_{2}^{\ast
}D_{2}=1$ and $Q_{2}=qB$ where $B$ is a free positive parameter. Since all
the braidings turn out to be trivial ( $\psi =\pi $) the coproducts,
antipodes and counits satisfy the ordinary Hopf algebra axioms. Choosing $B=%
\frac{1}{2}$\ the relation between $a$ and $a^{\ast }$ can equivalently be
written as

\begin{equation}
aa^{\ast }-a^{\ast }a=\cos \dfrac{(2N+1)\pi }{2k}
\end{equation}

\noindent and the spectrum is given by

\begin{equation}
a^{\ast }a\mid n\rangle =\dfrac{\sin \frac{n\pi }{k}}{2\sin \frac{\pi }{2k}}%
\mid n\rangle .
\end{equation}

\noindent This oscillator algebra, in the $k\rightarrow \infty $ limit,
reduces to the bosonic oscillator and hence we call it the generalized
exclusion algebra of the bosonic type. We note that the exchange phase is
constant (+1) and independent of $k$.

\subsection{Generalized exclusion of the fermionic type}

For the second solution we get $Q_{1}=-1$, and the exchange phase for the
creation and annihilation operators turn out to be -1 ( $z=-1$). The
oscillator algebra
\begin{eqnarray}  \label{eq:ferm2}
aa^{\ast }+a^{\ast }a=Q_{2}(q^{2N}-q^{-2(N+1)}),  \nonumber \\
aq^{N}=qq^{N}a,\quad q^{N}a^{\ast }=qa^{\ast }q^{N},  \nonumber \\
aq^{-N}=q^{-1}q^{-N}a,\quad q^{-N}a^{\ast }=q^{-1}a^{\ast }q^{-N}, \\
q^{N}q^{-N}=q^{-N}q^{N}=1_{A},\quad a^{k}=(a^{\ast })^{k}=0
\nonumber
\end{eqnarray}

\noindent has the Hopf algebra structure
\begin{eqnarray}
\Delta (q^{N})=D_{1}q^{N}\otimes q^{N},\ \ \Delta
(q^{-N})=D_{1}^{-1}q^{-N}\otimes q^{-N},  \nonumber \\
\Delta (a)=D_{2}a\otimes q^{N}+D_{1}^{-2}D_{2}q^{-N}\otimes a,  \nonumber \\
\Delta (a^{\ast })=D_{1}^{2}D_{2}^{\ast }a^{\ast }\otimes q^{N}+D_{2}^{\ast
}q^{-N}\otimes a^{\ast } \\
S(q^{N})=q^{-1}q^{-N},\quad S(q^{-N})=qq^{N},\quad S(a)=-q^{-1}a,\quad
S(a^{\ast })=-qa^{\ast },  \nonumber \\
\epsilon (q^{N})=D_{1}^{-1},\quad \epsilon (q^{-N})=D_{1},\quad \epsilon
(a)=\epsilon (a^{\ast })=0  \nonumber
\end{eqnarray}

\noindent where $D_{1}^{4}=q^{2}$, $D_{2}^{\ast }D_{2}=1$, $q=e^{i\frac{\pi
}{2k}}$and $Q_{2}=Bq^{1-k}=-Bq^{k+1}=-iqB$. The braiding relations are

\begin{equation}
\psi (a\otimes a)=-a\otimes a,\ \psi (a^{\ast }\otimes a^{\ast })=-a^{\ast
}\otimes a^{\ast },\ \psi (a\otimes a^{\ast })=-a^{\ast }\otimes a,\ \psi
(a^{\ast }\otimes a)=-a\otimes a^{\ast }.
\end{equation}

\noindent Identifying

\begin{equation}
a_{1}\equiv a\otimes 1_{A}, \ \ a_{2}\equiv 1_{A}\otimes a,\ \
a_{1}^{*}\equiv a^{*}\otimes 1_{A},\ \ a_{2}^{*}\equiv 1_{A}
\otimes a^{*}
\end{equation}

\noindent the braiding relations can equivalently be written as

\begin{equation}
a_{2}a_{1}=-a_{1}a_{2},\ \ a_{2}^{*}a_{1}^{*}=-a_{1}^{*}a_{2}^{*},\ \
a_{2}a_{1}^{*}=-a_{1}^{*}a_{2}, \ \ a_{2}^{*}a_{1}=-a_{1}a_{2}^{*}.
\end{equation}

\noindent The relation between $a$ and $a^{\ast }$ can equivalently be
written as

\begin{equation}
aa^{\ast }+a^{\ast }a=2B\sin \frac{(2N+1)\pi }{2k}
\end{equation}

\noindent and the spectrum is given by
\begin{equation}
a^{\ast }a\mid n\rangle =B\dfrac{\sin \frac{n\pi }{k}}{\cos \frac{\pi }{2k}}%
\mid n\rangle .
\end{equation}

\noindent Choosing $B=\frac{1}{\sqrt{2}}$ this oscillator gives the usual
fermion algebra for $k=2$ and hence we call it  the generalized exclusion
algebra of the fermionic type. We note that the exchange phase is constant
(-1) and independent of $k$.

\section{Covariance}

To study the covariance properties we start with the generalized oscillator
algebra (\ref{eq:osc1}) and make the transformation

\begin{equation}
\left(
\begin{array}{c}
a \\
a^{\ast } \\
q^{N} \\
q^{-N}
\end{array}
\right) ^{\prime }=\left(
\begin{array}{cccc}
t_{11} & t_{12} & t_{13} & t_{14} \\
t_{12}^{\ast } & t_{11}^{\ast } & t_{14}^{\ast } & t_{13}^{\ast } \\
t_{31} & t_{32} & t_{33} & t_{34} \\
t_{32}^{\ast }  & t_{31}^{\ast } & t_{34}^{\ast } & t_{33}^{\ast }
\end{array}
\right) \left(
\begin{array}{c}
a \\
a^{\ast } \\
q^{N} \\
q^{-N}
\end{array}
\right).
\end{equation}

\noindent The covariance of $\
(q^{N}q^{-N}=q^{-N}q^{N}=1_{A})^{\prime }$ gives

\be t_{31}=t_{32}= t_{34}=0,\quad t_{33}^{\ast }t_{33}=1 \ee

 \noindent and the covariance of $\
(aa^{\ast }-Q_{1}a^{\ast }a=Q_{2}q^{2N}+Q_{3}q^{-2N})^{\prime }$
gives
\begin{eqnarray}  \label{eq:cov1}
t_{11}t_{12}^{\ast }-Q_{1}t_{12}^{\ast }t_{11}&=&0,  \nonumber \\
t_{12}t_{11}^{\ast }-Q_{1}t_{11}^{\ast }t_{12}&=&0,  \nonumber \\
Q_{2}t_{11}t_{11}^{\ast }-Q_{1}Q_{2}t_{12}^{\ast }t_{12}+t_{13}t_{14}^{\ast
}-Q_{1}t_{14}^{\ast }t_{13}&=&Q_{2}t_{33}^{2},  \nonumber \\
Q_{3}t_{11}t_{11}^{\ast }-Q_{1}Q_{3}t_{12}^{\ast }t_{12}+t_{14}t_{13}^{\ast
}-Q_{1}t_{13}^{\ast }t_{14}&=&Q_{3}(t_{33}^{\ast })^{2},  \nonumber \\
t_{11}t_{11}^{\ast }-t_{11}^{\ast }t_{11}+Q_{1}^{-1}t_{12}t_{12}^{\ast
}-Q_{1}t_{12}^{\ast }t_{12}&=&0, \\
qt_{13}t_{11}^{\ast }-Q_{1}t_{11}^{\ast }t_{13}+t_{12}t_{14}^{\ast
}-qQ_{1}t_{14}^{\ast }t_{12}&=&0,  \nonumber \\
q^{-1}t_{11}t_{13}^{\ast }-Q_{1}t_{13}^{\ast }t_{11}+t_{14}t_{12}^{\ast
}-q^{-1}Q_{1}t_{12}^{\ast }t_{14}&=&0,  \nonumber \\
qt_{11}t_{14}^{\ast }-Q_{1}t_{14}^{\ast }t_{11}+t_{13}t_{12}^{\ast
}-qQ_{1}t_{12}^{\ast }t_{13}&=&0,  \nonumber \\
q^{-1}t_{14}t_{11}^{\ast }-Q_{1}t_{11}^{\ast }t_{14}+t_{12}t_{13}^{\ast
}-q^{-1}Q_{1}t_{13}^{\ast }t_{12}&=&0,  \nonumber \\
t_{13}t_{13}^{\ast }+t_{14}t_{14}^{\ast }-Q_{1}t_{13}^{\ast
}t_{13}-Q_{1}t_{14}^{\ast }t_{14}&=&0.  \nonumber
\end{eqnarray}

\noindent The compatibility of these relations with the *-operation gives $%
Q_{1}=\pm 1$ and $Q_{3}=Q_{2}^{\ast }$. For $Q_{1}=q=1$ the entries of the
transformation matrix commute among themselves and setting $t_{33}=1$ and $%
t_{13}=t_{14}=0$ this gives the group $SU(1,1)$. For $q$ a root of unity the
covariance of the nilpotency condition $(a^{\prime })^{k}=0$ gives $t_{12}=0$%
. For $Q_{1}=-1$ the last equation in (\ref{eq:cov1}) requires $%
t_{13}=t_{14}=0$ and the transformation for the generalized exclusion of the
fermionic type given by (\ref{eq:ferm2}) turns out to be trivial. The
algebra for the generalized exclusion of the bosonic type given by (\ref
{eq:bos1}) is covariant under the transformation

\begin{equation}
\left(
\begin{array}{c}
a \\
a^{\ast } \\
q^{N} \\
q^{-N}
\end{array}
\right) ^{\prime }=\left(
\begin{array}{cccc}
t_{11} & 0 & t_{13} & t_{14} \\
0 & t_{11}^{\ast } & t_{14}^{\ast } & t_{13}^{\ast } \\
0 & 0 & t_{33} & 0 \\
0 & 0 & 0 & t_{33}^{\ast }
\end{array}
\right) \left(
\begin{array}{c}
a \\
a^{\ast } \\
q^{N} \\
q^{-N}
\end{array}
\right)
\end{equation}

\noindent if the entries of the transformation matrix satisfy
\begin{eqnarray}
t_{11}t_{11}^{\ast }&=&t_{11}^{\ast }t_{11},  \nonumber \\
t_{11}t_{13}&=&q^{3}t_{13}t_{11},  \nonumber \\
t_{11}t_{13}^{\ast }&=&qt_{13}^{\ast }t_{11},  \nonumber \\
t_{11}t_{14}&=&q^{-3}t_{14}t_{11},  \nonumber \\
t_{11}t_{14}^{\ast }&=&q^{-1}t_{14}^{\ast }t_{11},  \nonumber \\
t_{11}t_{33}&=&t_{33}t_{11},  \nonumber \\
t_{11}t_{33}^{\ast }&=&t_{33}^{\ast }t_{11},  \nonumber \\
t_{13}t_{14}&=&q^{-4}t_{14}t_{13},  \nonumber \\
t_{13}t_{14}^{\ast }-t_{14}^{\ast }t_{13}+Q_{2}t_{11}t_{11}^{\ast
}&=&Q_{2}t_{33}^{2}, \\
t_{13}t_{13}^{\ast }-t_{14}^{\ast }t_{14}&=&t_{13}^{\ast
}t_{13}-t_{14}t_{14}^{\ast },  \nonumber \\
t_{13}t_{33}&=&qt_{33}t_{13},  \nonumber \\
t_{13}t_{33}^{\ast }&=&q^{-1}t_{33}^{\ast }t_{13},  \nonumber \\
t_{14}t_{33}&=&qt_{33}t_{14},  \nonumber \\
t_{14}t_{33}^{\ast }&=&q^{-1}t_{33}^{\ast }t_{14},  \nonumber \\
t_{33}^{\ast }t_{33}&=&1,  \nonumber \\
t_{13}^{k}&=&t_{14}^{k}=0  \nonumber
\end{eqnarray}

\noindent and $*$-conjugates. The deformation parameter $q$ is a root of
unity, $q=e^{i\frac{\pi }{2k}}$. The coproducts
\begin{eqnarray}
\Delta (t_{11})=t_{11}\otimes t_{11},\ \ \Delta (t_{13})=t_{11}\otimes
t_{13}+t_{13}\otimes t_{33},  \nonumber \\
\Delta (t_{33})=t_{33}\otimes t_{33},\ \ \Delta (t_{14})=t_{11}\otimes
t_{14}+t_{14}\otimes t_{33}^{\ast },  \nonumber \\
\Delta (t_{11}^{*})=t_{11}^{*}\otimes t_{11}^{*},\ \ \Delta
(t_{13}^{*})=t_{11}^{*}\otimes t_{13}^{*}+t_{13}^{*}\otimes t_{33}^{*}, \\
\Delta (t_{33}^{*})=t_{33}^{*}\otimes t_{33}^{*},\ \ \Delta
(t_{14}^{*})=t_{11}^{*}\otimes t_{14}^{*}+t_{14}^{*}\otimes t_{33},
\nonumber
\end{eqnarray}
\noindent the antipodes
\begin{eqnarray}
S(t_{11})=t_{11}^{-1}\quad S(t_{13})=-t_{11}^{-1}t_{13}t_{33}^{-1}\quad
S(t_{14})=-t_{11}^{-1}t_{14}t_{33}\quad S(t_{33})=t_{33}^{-1}, \\
S(t_{11}^{\ast })=(t_{11}^{\ast })^{-1},\quad S(t_{13}^{\ast
})=-(t_{11}^{\ast })^{-1}t_{13}^{\ast }t_{33},\quad S(t_{14}^{\ast
})=-(t_{11}^{\ast })^{-1}t_{14}^{\ast }t_{33}^{-1},\quad S(t_{33}^{\ast
})=t_{33}  \nonumber
\end{eqnarray}
\noindent and counits

\begin{eqnarray}
\epsilon (t_{11})=\epsilon (t_{11}^{\ast })=\epsilon (t_{33})=\epsilon
(t_{33}^{\ast })=1, \\
\epsilon (t_{13})=\epsilon (t_{13}^{\ast })=\epsilon (t_{14})=\epsilon
(t_{14}^{\ast })=0  \nonumber
\end{eqnarray}

\noindent complete the Hopf algebra structure.

\section{Conclusion}

In this work we investigated the braided Hopf algebra structure of
generalized oscillator with generalized exclusion $a^{k}=0$, i.e.,
there can be at most $k-1$ particles at the same state. We find
that there are two solutions one of which is the generalized
exclusion of the bosonic type where the exchange phase is $(+1)$
and the limit $k\rightarrow \infty $ gives the bosonic oscillator.
The second solution is the generalized exclusion of the fermionic
type where the exchange phase is $(-1)$ and the limit
$k\rightarrow 2$ gives the fermionic oscillator. We note that in
the operator algebra with exclusion $a^{k}=0$ we start with the
generalized
exchange phases (braidings) but the Hopf algebra structure chooses only the $%
\pm 1$ factor. If we do not impose any exclusion on the
creation-annihilation operators it is possible to find Hopf algebra
solutions with the generalized exchange phase (braidings) \cite{art4}.
Hence, in the context of Hopf algebras, we arrive at the conclusion that the
generalization of exclusion does not necessarily imply the generalization of
the exchange phase and vice versa.

The Hopf algebras we found for the oscillators can be used to generate the
Hopf structure of Lie algebras through a realization using some extension of
the Jordan-Schwinger map just like the ones discussed for parabosonic and
parafermionic algebras \cite{art14}.

\end{document}